\documentclass[sn-mathphys-num]{sn-jnl}

\usepackage{graphicx}%
\usepackage{multirow}%
\usepackage{amsmath,amssymb,amsfonts}%
\usepackage{amsthm}%
\usepackage{mathrsfs}%
\usepackage[title]{appendix}%
\usepackage{xcolor}%
\usepackage{textcomp}%
\usepackage{manyfoot}%
\usepackage{booktabs}%
\usepackage{algorithm}%
\usepackage{algorithmicx}%
\usepackage{algpseudocode}%
\usepackage{listings}%

\usepackage[utf8]{inputenc}

\usepackage{hyperref}
\usepackage{cleveref}
\usepackage{array}
\usepackage{makecell}
\usepackage{mathtools}
\usepackage{thmtools}
\usepackage{thm-restate}
\usepackage{comment}
\usepackage{soul}
\usepackage{tikz}
\usepackage{adjustbox}
\usetikzlibrary{positioning}

\usepackage{caption}
\usepackage{float}

\numberwithin{equation}{section}

\theoremstyle{definition}
\newtheorem{theorem}{Theorem}[section]
\newtheorem*{theorem*}{Theorem}
\newtheorem{example}[theorem]{Example}
\newtheorem*{example*}{Example}
\newtheorem{lemma}[theorem]{Lemma}
\newtheorem*{lemma*}{Lemma}
\newtheorem{corollary}[theorem]{Corollary}
\newtheorem*{corollary*}{Corollary}

\newtheorem*{definition*}{Definition}
\newtheorem{proposition}[theorem]{Proposition}
\newtheorem*{proposition*}{Proposition}
\newtheorem{remark}[theorem]{Remark}
\newtheorem*{remark*}{Remark}
\newtheorem{conjecture}[theorem]{Conjecture}

\usepackage{ mathrsfs }

\newenvironment{amatrix}[1]{%
  \left(\begin{array}{@{}*{#1}{c}|c@{}}
}{%
  \end{array}\right)
}

\begin{document}

\title[Periods and Atomic Firing Sequences of Parallel Chip-Firing Games on Directed Graphs]{Periods and Atomic Firing Sequences of Parallel Chip-Firing Games on Directed Graphs}


\author[1]{\fnm{David} \sur{Ji}}\email{david.ji@mtsdstudent.us}
\equalcont{These authors contributed equally to this work.}

\author[2]{\fnm{Michael} \sur{Li}}\email{michaelli20180101@gmail.com}
\equalcont{These authors contributed equally to this work.}

\author[3]{\fnm{Daniel} \sur{Wang}}\email{danielw25@lakesideschool.org}
\equalcont{These authors contributed equally to this work.}

\affil*[1]{\orgname{Montgomery High School}, \orgaddress{\city{Skillman}, \postcode{08558}, \state{New Jersey}, \country{United States of America}}}

\affil[2]{\orgname{University High School}, \orgaddress{\city{Irvine}, \postcode{92612}, \state{California}, \country{United States of America}}}

\affil[3]{\orgname{Lakeside School}, \orgaddress{\city{Seattle}, \postcode{98125}, \state{Washington}, \country{United States of America}}}


\abstract{In 1992, Bitar and Goles introduced the parallel chip-firing game on undirected graphs. Two years later, Prisner extended the game to directed graphs. While the properties of parallel chip-firing games on undirected graphs have been extensively studied, their analogs for parallel chip-firing games on directed graphs have been sporadic. In this paper, we prove the outstanding analogs of the core results of parallel chip-firing games on undirected graphs for those on directed graphs. We find the possible periods of a parallel chip-firing game on a directed simple cycle and introduce the method of Gauss-Jordan elimination on a Laplacian-like matrix to establish a lower bound on the maximum period of a parallel chip-firing game on an orientation of an undirected complete graph and an undirected complete bipartite graph. Finally, we expand the method of motors by Jiang, Scully, and Zhang to directed graphs to show that a binary string $s$ can be the atomic firing sequence of a vertex in a parallel chip-firing game on a strongly connected directed graph if and only if $s$ contains $1$ or $s=0$.}




\maketitle

\subsection*{Statements and Declarations}
The authors declare that they have no known competing financial interests or personal relationships that could have appeared to influence the work reported in this paper.

\section{Introduction}\label{intro}

 The \emph{parallel chip-firing game} is an automaton on a finite, connected graph $G(V,\,E)$  that was first introduced by Bitar and Goles \cite{bitarandgoles} in 1992. A parallel chip-firing game begins with some \emph{chips} distributed among each $v\in V$. On each subsequent round, all vertices with at least as many chips as neighbors simultaneously \emph{fire} one chip to each of their neighbors. All other vertices \emph{wait}. As observed by Bitar and Goles \cite{bitarandgoles}, all parallel chip-firing games are periodic with some minimal period length $T$, because the quantity of chips and the number of vertices in any parallel chip-firing game are fixed.

The parallel chip-firing game on undirected graphs has been extensively studied for trees \cite{bitarandgoles}, simple cycles \cite{dallasta}, complete graphs \cite{levine}, and complete bipartite graphs \cite{jiang}. See Table \ref{periods} for the possible period lengths of parallel chip-firing games on these graphs. For general graphs, Kominers and Kominers \cite{kominers} showed that $T=1$ for all parallel chip-firing games with enough chips. Later, Bu, Choi, and Xu \cite{yunseo} made this range of chips exact. Furthermore, Kiwi et al. \cite{kiwi} showed that $T$ cannot be bounded by a polynomial in $|V|$. Lastly, Scully, Jiang, and Zhang \cite{scully} distinguished the possible \emph{atomic firing sequences} --- length-$T$ binary strings that indicate the rounds in which a vertex fires within a periodic sequence --- of a vertex in a parallel chip-firing game. 

In 1994, Prisner \cite{prisner} extended the parallel chip-firing game to directed graphs. Prisner \cite{prisner} found that $T=1$ for all games on acyclic directed graphs and that there exists no polynomial bound on $|V|$ for the period of a parallel chip-firing game on a directed graph. Later, Goles and Prisner \cite{sourcereversal} found that $T\leq |V|$ for all parallel chip-firing games on orientations of $\overline{K_n}*H$, where $\overline{K_n}$ is an empty graph, $*$ denotes the graph join operation, and $H$ is an arbitrary graph. Finally, Ndoundam, Tchuente, and Tadonki \cite{hypercube} characterized the possible periods on orientations of the $n$-cube graph.

In this paper, we prove the analogs of results for the parallel chip-firing game on undirected graphs for those on directed graphs. In Section \ref{prelim}, we establish preliminaries. In Section \ref{parallel}, we characterize the possible periods of parallel chip-firing games on certain classes of directed graphs. In addition, we observe that while there always exists a parallel chip-firing game with $T=2$ on any undirected graph, the same is not true for orientations of complete directed graphs with at most four vertices. In Section \ref{atomic}, we prove that a binary string $s$ can be the atomic firing sequence of a vertex in a parallel chip-firing game on a strongly connected directed graph if and only if $s$ contains 1 or $s=0$. Finally, in Section \ref{futuredirections}, we present unsolved conjectures and areas for future study.

\begin{table}
    \centering
    \begin{tabular}{|l|m{10em}|m{15em}|}
    \hline
        Graph & \makecell{Undirected} & \makecell{Directed} \\
         \hline
         \hline
         Path $P_n$ & $T=1,\,2$; \cite{bitarandgoles} & $T=1$; \cite{prisner} \\
         \hline
         Non-path tree & $T=1,\,2$; \cite{bitarandgoles} & $T=1$; \cite{prisner} \\
         \hline
         Simple cycle $C_n$ & $T=2$ or $T=i$ for all $i|n$; \cite{dallasta} & $T=i$ for all $i|n$; \newline \Cref{cycles} \\
         \hline
          Complete $K_n$ & $T=i$ for all $i\leq n$; \cite{levine} & Maximum possible $T$ grows with at least $(n-1)!$; \Cref{tournaments} \\
         \hline
         Complete bipartite $K_{a,\,b}$ & $T=i$ or $2i$ for all $i\leq \min(a,\,b)$; \cite{jiang} & Maximum possible $T$ grows with at least $\min(a,\,b)!$; \Cref{bipartites} \\
         \hline
    \end{tabular}
    \caption{The period lengths $T$ of parallel chip-firing games on special classes of graphs.}
    \label{periods}
\end{table}

\section{Preliminaries}\label{prelim}

\begin{figure}[H]
    \centering
    \begin{tikzpicture}
      [every node/.style={circle, draw, minimum size=6mm, inner sep=1mm, font=\footnotesize},
       line cap=round, line join=round, thick]
    
      \node (1) at (0, 0) {1};
      \node (2) at (2, 0) {2};
      \node (3) at (1, -1.5) {3};
      \node (4) at (3, -1.5) {4};
    
      \draw[-] (1) -- (2);
      \draw[-] (2) -- (3);
      \draw[-] (3) -- (4);
      \draw[-] (1) -- (3);
    \end{tikzpicture}
    \hspace{2cm}
    \begin{tikzpicture}
      [every node/.style={circle, draw, minimum size=6mm, inner sep=1mm, font=\footnotesize},
       ->, line cap=round, line join=round, thick]
    
      \node (1) at (0, 0) {1};
      \node (2) at (2, 0) {2};
      \node (3) at (1, -1.5) {3};
      \node (4) at (3, -1.5) {4};
    
      \draw[->] (1) -- (2);
      \draw[->] (2) -- (3);
      \draw[->] (3) -- (1);
      \draw[->] (3) -- (4);
    \end{tikzpicture}
    \caption{An undirected graph $G$ (left) and an orientation of $G$.}
    \label{orientation}
\end{figure}
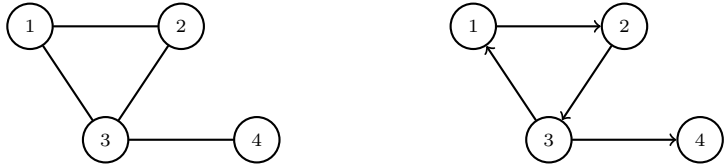

We define a \emph{directed graph} $D=(V,\,E)$ as a set of vertices $V$ and a set of ordered pairs of vertices (\emph{edges}) $E$. The edge $(u,\,v)$ for $u,\,v\in V$ is said to be \emph{directed from $u$ to $v$}. We define an \emph{orientation} of an undirected graph $G=(V,\,E)$ to be a directed graph $D=(V,\,E')$ such that either $(u,\,v)\in E'$ or $(v,\,u)\in E'$ if and only if there exists an edge between $u$ and $v$ in $G$, but there do not exist vertices $u$ and $v$ such that $(u,\,v)$ and $(v,\,u)$ are both elements in $E'$. The \emph{in-degree} of a vertex $v$ is denoted $\deg^-v$ and is the number of edges directed from another vertex to $v$. The \emph{out-degree} of a vertex $v$ is denoted $\deg^+v$ and is the number of edges directed from $v$ to another vertex. In the directed graph in Figure \ref{orientation}, $\deg^+v_3=2$ and $\deg^-v_3=1$. The \emph{degree} of a vertex $v$ is denoted $\deg v$ and is the total number of edges including $v$. In Figure \ref{orientation}, $\deg v_3=3$. The \emph{distance} between two vertices $u$ and $v$ is the length of the shortest path from $u$ to $v$ and is denoted $d(u,\,v)$. In the undirected graph in Figure \ref{orientation}, $d(v_1,\,v_4)=2$. In the directed graph, $d(v_1,\,v_4)=3$. For parallel chip-firing games on directed graphs, a vertex fires when the number of chips it possesses exceeds its out-degree, upon which the vertex fires one chip across each outgoing edge.

Bitar and Goles \cite{bitarandgoles} noted that the fixed numbers of chips and vertices make every parallel chip-firing game on a finite undirected graph periodic. Prisner \cite{prisner} noted that the same is true for parallel chip-firing games on finite directed graphs. We say that a parallel chip-firing game on a directed graph has period $T$ if $T$ is the smallest integer for which there exists a round $t_0$ such that every vertex has the same number of chips in round $t$ as in round $t+T$ for all $t \geq t_0$. We denote the least such $t_0$ as round 0. Let $c_t(v)$ be the number of chips on a vertex $v$ at the beginning of round $t$. Let $F_t(v)$ be $1$ if $v$ fires on round $t$ and $0$ if it does not. Let $f_v=\sum_{t=0}^{T-1} F_t(v)$ be the total number of times $v$ fires in the first $T$ rounds. 

\begin{figure}
    \centering
    \trimbox{0cm 1cm 0cm 0cm}{
    \begin{tikzpicture}
      [every node/.style={circle, draw, minimum size=6mm, inner sep=0},
       scc/.style={draw=black, thick, rounded corners, inner sep=2mm}]
    
      \node (A) at (-6, 1) {A};
      \node (B) at (-5, 1) {B};
      \node (C) at (-5, 0) {C};
      \node (D) at (-6, 0) {D};
      \node (E) at (-2.5, 0) {E};
      \node (F) at (-2, 0.7) {F};
      \node (G) at (-1.5, 0) {G};
      \node (H) at (-3.7, 1.5) {H}; 
      
      \draw[->] (A) -- (B);
      \draw[->] (B) -- (C);
      \draw[->] (C) -- (D);
      \draw[->] (D) -- (A);
      \draw[->] (C) -- (E);
      \draw[->] (E) -- (F);
      \draw[->] (F) -- (G);;
      \draw[->] (G) -- (E);
      \draw[->] (B) -- (H);
      \draw[->] (F) -- (H);
    
      \begin{scope}[]
        \draw[scc] (-6.5, -0.5) rectangle (-4.5, 1.5);
        \draw[scc] (-3, -0.5) rectangle (-1, 1.2);
        \draw[scc] (-4.2, 1) rectangle (-3.2, 2); 
      \end{scope}
    
      \node (SCC1) at (3, 0) {SCC1};
      \node (SCC2) at (5, 0) {SCC2};
      \node (SCC3) at (4, 1) {SCC3};
      
      \draw[->] (SCC1) -- (SCC2);
      \draw[->] (SCC1) -- (SCC3);
      \draw[->] (SCC2) -- (SCC3);
    
      \node[draw=none] at (-4, -1) {Original Graph};
      \node[draw=none] at (4, -1) {Condensation};
    
    \end{tikzpicture}
    }
    \caption{A directed graph (left) with strongly connected components boxed and its condensation (right).}
    \label{condensation}
\end{figure}
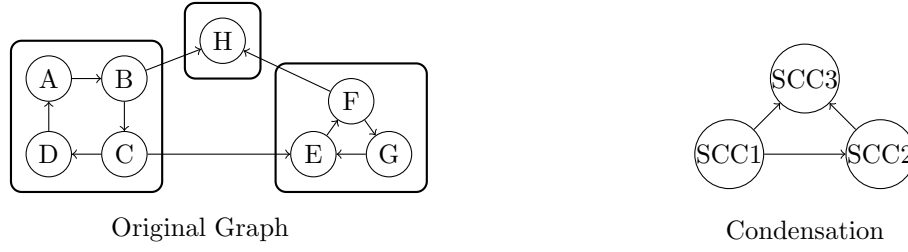

We say that a directed graph $D=(V,\,E)$ is \emph{strongly connected} if there exists a path in both directions between every pair of vertices $v,\,u\in V$. A \emph{strongly connected component} of $D$ is a subgraph of $D$ that is strongly connected and maximal: no more vertices or edges may be included while keeping the subgraph strongly connected. We let the \emph{condensation} of a directed graph $D$ be a directed graph with all strongly connected components of $D$ as vertices. Borrowing notation from Prisner \cite{prisner}, we define a \emph{sink component} as a strongly connected component that forms a sink in the condensation of $D$. In Figure \ref{condensation}, SCC3 is a sink component. We say that a vertex is \emph{forever passive} if it does not fire in all future states.

Prisner \cite{prisner} found that only sink components are relevant to the behavior of parallel chip-firing games on a directed graph.

\begin{proposition} [{\cite[Proposition 2.2]{prisner}}]
\label{nonsink}
     In every round $t\geq 0$ in a parallel chip-firing game on $D$, all vertices in nonsink components are forever passive.
\end{proposition}

By \Cref{nonsink}, in any parallel chip-firing game on a directed graph with multiple sink components, the number of chips in each sink component is eventually fixed. The game can then be considered a superposition of independent parallel chip-firing games on each sink component. In this paper, we will therefore assume that all graphs are strongly connected (i.e. the only sink component of a graph $D$ is $D$ itself) unless otherwise stated. We first prove that no vertices are forever passive in any parallel chip-firing game on a strongly connected directed graph.

\begin{lemma}
\label{lost}
    In any parallel chip-firing game on a strongly connected directed graph $D$ with $T>1$, no vertices are forever passive.
\end{lemma}

\begin{proof}
    As $T>1$, at least one vertex $v'$ must fire in each period. Since $c_0(v)=c_T(v)$ for all $v\in V$, each of the vertices which $v'$ fires to must also fire in each period. Similarly, all the vertices that are distance 2 from $v'$ must fire in each period. In general, every vertex $u$ for which there exists a directed path from $v'$ to $u$ must fire in each period. Since $D$ is strongly connected, such a path exists to any vertex $u\in V$. Thus, all vertices fire at least once in each period.
\end{proof}

\section{Period Lengths of Parallel Chip-Firing Games}\label{parallel}

\subsection{Trees}

The possible periods of parallel chip-firing games on directed trees were identified by Prisner \cite{prisner}. We begin by citing his result.

\begin{proposition} [{\cite[Corollary 2.3]{prisner}}]
\label{acyclic}
    $T=1$ for all parallel chip-firing games on directed acyclic graphs.
\end{proposition}

For completeness, we state the following direct corollaries of \Cref{acyclic}.

\begin{corollary}
\label{path_par}
      $T=1$ for all parallel chip-firing games on directed paths.
\end{corollary} 

\begin{corollary}
\label{non-path_par}
    $T=1$ for all parallel chip-firing games on directed non-path trees.
\end{corollary}

\subsection{Cycles}
Consider the simple directed cycle $C_n$. We label the vertices of the cycle $v_0, v_1, v_2,\dots, v_{n-1}$ such that there exists an edge directed from $v_{i+1}$ to $v_i$ for all $0 \leq i \leq n-1$, with indices taken modulo $n.$ We first show that if any vertex ever has more than 1 chip, the period of the game is 1.

\begin{lemma}
    \label{c=1}
    In a parallel chip-firing game on $C_n$, if $c_t(v)>1$ for any $v$, then $T=1$.
\end{lemma}

\begin{proof}
    Suppose $c_t(v_0)>1$. Note that if $c_t(v)\geq1$ then $c_{t+1}(v) \leq c_t(v)$, so we must have $1 < c_t(v_0)=c_{t+1}(v_0)=\dots=c_{t+T-1}(v_0)=c_{t+T}(v_0)$. We conclude that $v_0$ must fire on every round and so $v_1$ must fire on every round as well. By induction, we conclude that every vertex fires on every round, so $T=1$.
\end{proof}

We now prove our main theorem about directed cycles.

\begin{theorem}
    \label{cycles}
    The possible periods $T$ of a parallel chip-firing game on $C_n$ are $T=i$ for all $i$ such that $i \; \vert \; n.$
\end{theorem}

\begin{proof}
    Consider a parallel chip-firing game on $C_n$ with period $T$. If $T=1$, we have $T\;\vert\; n$. Otherwise, by \Cref{c=1}, we have $c_t(v)\leq 1$ for all $v$ and $t$. We conclude that $c_t(v)=F_t(v)$ for all $v$ and $t$.
    
    Note that if $F_{t}(v_i)=1,$ then $v_{i-1}$ has one chip after round $t,$ so $F_{t+1}(v_{i-1})=1.$ Thus, if $F_t(v)=1,$ then $F_{t+n}(v)=1.$ Additionally, if $F_{t+n}(v)=1,$ then $F_{t+nT}(v)=1$, and by the definition of a period, $F_t(v)=1$ as well. Therefore, $F_t(v)=F_{t+n}(v)$ for all $v$ and $t$, so $c_t(v)=c_{t+n}(v)$ for all $v$ and $t$. Thus, we must have $T \; \vert \; n.$

    To construct a parallel chip-firing game with period $T=i$ for any $i \; \vert \; n,$ place a single chip on each vertex $v_j$ that satisfies $j \equiv 0 \pmod{i}$. On round $t,$ there will be a chip on each vertex $v_j$ that satisfies $j \equiv -t \pmod{i}.$ The smallest $T$ for which $c_v(t) = c_v(t+T)$ is therefore $T = i$.
\end{proof}

\subsection{Complete Graphs}

We first note as a surprising aside that while there exists a parallel chip-firing game on any undirected graph with $T=2$, it is possible for there to be no orientations of a graph that allow for a parallel chip-firing game with $T=2$. We begin by citing a lemma from Jiang \cite{jiang}.

\begin{lemma}[{\cite[Proposition 2.5]{jiang}}] 
\label{constance}
    For any parallel chip-firing game on an undirected graph, $c_T(v)=c_0(v)$ for all $v\in V$ if and only if $\sum_{t=0}^{T-1} F_t(v)$ is constant for all $v\in V$.
\end{lemma}

In \Cref{talways2}, we show that there exists a parallel chip-firing game on any undirected graph $G$ such that every vertex fires exactly once in the first 2 rounds. However, parallel chip-firing games on directed graphs behave differently, as demonstrated by the following lemma.

\begin{lemma}
\label{tnot2}
    No parallel chip-firing games on any orientation of $K_n$ with $n\leq 4$ have $T=2$.
\end{lemma}

\begin{proof}
    Note that $K_2$ is a path, and $K_3$ is a cycle. For these cases, see \Cref{path_par} and \Cref{cycles}, respectively. Then consider a parallel chip-firing game on a strongly connected orientation $D=(V,\,E)$ of $K_4$ with $T=2$. Recall that because $D$ is strongly connected, each vertex $v$ in $D$ must satisfy $\deg^+v>0$ and $\deg^-v>0$. Also, by Lemma \ref{lost}, each vertex must fire at least once in each period.
    
    Since $\deg v=3$, $v$ must be one of two types: $\deg^- v =1$ and $\deg^+v=2$, or $\deg^-v=2$ and $\deg^+v=1$. Additionally, since $\sum_{v\in V} \deg^-v = \sum_{v\in V} \deg^+v$, there must be two vertices of each type. Let $\deg^+ v_1= \deg^+ v_2=1$, and let $\deg^+ u_1=\deg^+ u_2=2$.
    
    Note that there must exist an edge between $v_1$ and $v_2$ and that regardless of its direction, one of the two vertices must view it as its single out-edge. Thus, there is only one edge directed from $v_1$ or $v_2$ to $u_1$ or $u_2$. Without loss of generality, let it be directed from $v_1$ to $u_1$. On the other hand, there are three edges directed from $u_1$ or $u_2$ to $v_1$ or $v_2$.
    
    Since $T=2$, we must have $f_{v_1}\leq 2$. However, $f_{u_1}\geq 1$ and $f_{u_2}\geq 1$. Then there must be at most two chips that travel from a vertex in $V=\{v_1,\,v_2\}$ to one in $U=\{u_1,\,u_2\}$, but at least three chips that travel from a vertex in $U$ to one in $V$. It is impossible for every vertex to have the same number of chips after 2 rounds, so $T\neq 2$ for all parallel chip-firing games on $D$.
\end{proof}

Next, let $\tau(n)$ be the maximum possible period of a parallel chip-firing game on an orientation of $K_n$. We show that $\tau(n)$ is at least asymptotic to $(n-1)!$. We conjecture that $\tau(n)$ is indeed asymptotic to $(n-1)!$

\begin{conjecture}
    \label{conj_tourney}
    It holds that $\tau(n)= \Theta((n-1)!)$.
\end{conjecture}

We define the following orientation of $K_n$ with vertices $v_1,\,v_2,\, \ldots,\, v_n$ to be a \emph{useful} orientation.
\begin{itemize}
    \item There exist edges directed from $v_i$ to $v_{i+1}$ for each $1 \leq i < n$, and there exists an edge directed from $v_n$ to $v_1$.
    \item There exist edges directed from $v_1$ to every other vertex except $v_n$.
    \item There exist edges directed from $v_j$ to all vertices $v_k$ with $1 < k \leq j-2$.
\end{itemize}

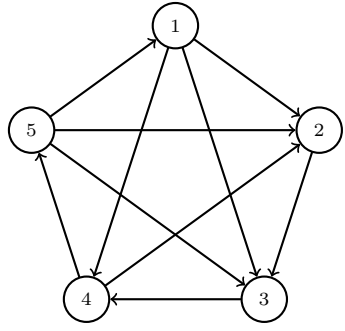
\begin{figure}
    \centering
    \begin{tikzpicture}
      [every node/.style={circle, draw, minimum size=6mm, inner sep=1mm, font=\footnotesize},
       ->, line cap=round, line join=round, thick]
    
      \node (1) at (90:2) {1};
      \node (5) at (162:2) {5};
      \node (4) at (234:2) {4}; 
      \node (3) at (306:2) {3};
      \node (2) at (18:2) {2}; 
    
      \draw[->] (1) -- (2);
      \draw[<-] (5) -- (4);
      \draw[<-] (4) -- (3);
      \draw[<-] (3) -- (2);
      \draw[->] (1) -- (4);
      \draw[->] (1) -- (3);
      \draw[<-] (3) -- (5);
      \draw[<-] (2) -- (5);
      \draw[<-] (2) -- (4);
      \draw[->] (5) -- (1);  
    \end{tikzpicture}
    \caption{A useful orientation of $K_5$.}
    \label{tournament}
\end{figure}

See Figure \ref{tournament} for a diagram of a useful orientation of $K_5$. Unless otherwise stated, all future references of $K_n$ will refer specifically to the useful orientation of $K_n$.

The following lemma shows that vertex 3 fires the most in all parallel chip-firing games on $K_n$ with $n\geq 4$.

\begin{lemma}
    \label{vertex3}
    For any $n \geq 4,$ all parallel chip-firing games on $K_n$ satisfy $f_3\geq f_i$ for all $1\leq i\leq n$.
\end{lemma}

\begin{proof}
     Since $\deg^+ v_3=1$, it must fire once in each period for every chip received in that period. Note that $v_3$ receives chips from all other vertices except $v_4$. Thus, $f_3\geq f_i$ for all $i\neq 4$. Also note that $\deg^+ v_2=1$, and that $v_2$ receives chips from $v_4$. Then $f_2\geq f_4$. Since $f_3\geq f_2$, we have $f_3\geq f_4$.
\end{proof}

Let the \emph{$c$-convergent period} $T_c(G)$ of a graph $G$ be $\min(T)$ over all parallel chip-firing games on $G$ with at least $c$ chips. Note that $T_c(G)=1$ for any undirected graph $G$, regardless of $c$ \cite{yunseo}. We will show that the same is not true for directed graphs.

\begin{theorem}
    \label{tournaments}
    For the useful orientation of $K_n$, with $c=\frac{n(n-3)}{2}+1$, we have $T_c(K_n)=\Omega((n-1)!)$.
\end{theorem}

\begin{proof}
    Since the number of chips on each vertex does not change after a full period, the following equations must hold.
    \begin{equation} \label{complete_original_system}
    \begin{aligned}
        & f_n = (n-2)f_1 \\
        & f_1+f_4+ \ldots + f_n = f_2 \\
        &f_1+f_{j-1} + f_{j+2} + f_{j+3} + \ldots + f_n=(j-2)f_j \quad\text{for}\quad 3 \leq j \leq n-2\\
        &f_1+f_{n-2}=(n-3)f_{n-1}\\
        &f_{n-1}=(n-2)f_n
    \end{aligned}
    \end{equation}
    
    For any solution, scaling gives a new solution, so System \ref{complete_original_system} is not independent. We arbitrarily drop the first equation, leaving a system of $n-1$ independent equations in $n$ variables. These remaining equations are linearly independent, as can be seen later in the un-augmented form of Matrix \ref{complete_simplified_matrix}, whose rank is equal to the number of rows.
    
    By definition, each of $f_1, f_2, \ldots,\, f_n \geq 0.$ By \Cref{lost}, if any $f_i=0$, we have $T=1$ and all vertices are forever passive. However, the largest possible $c$ for which all vertices can be forever passive is $\sum_v (\deg^+ v-1)=|E|-|V|=\frac{n(n-3)}{2}<c$. Thus, we further constrain $f_i>0$ for all $i$. Note that for any solution $(f_1, f_2, \ldots ,\,f_n)$, it must be true that $(\lambda f_1, \lambda f_2, \ldots,\, \lambda f_n)$ is also a solution. To find the minimal solution, we set $f_1=1.$ System \ref{complete_original_system} has now been reduced to the following.
    \begin{equation} \label{complete_simplified_system}
    \begin{aligned}
        &f_2-f_4-\ldots-f_n=1\\
        &-f_{j-1}+(j-2)f_j - f_{j+2} - f_{j+3} - \ldots - f_n=1 \quad\text{for}\quad j: \; 3 \leq j \leq n-2 \\
        &-f_{n-2}+(n-3)f_{n-1}=1\\
        &-f_{n-1}+(n-2)f_n=0
    \end{aligned}
    \end{equation}
     
    Note that these equations are linearly independent. Therefore, we can let a well-defined $T_n$ denote the maximal value of $f_i$ satisfying System \ref{complete_simplified_system}. By \Cref{vertex3}, for any  solution with $n\geq 4$, we have $f_3\geq f_i$ for all $i$. Therefore, $T_n=f_3$ if $n\geq 4$. We claim that, for $n\geq 4$, $$T_n=T_{n-2}+(n-1)T_{n-1}.$$
     
    System \ref{complete_simplified_system} can be written as the following augmented matrix $M_n=(A|B)$, where $A$ is an $(n-1)\times (n-1)$ matrix and $B$ is a $(n-1)\times 1$ column vector. We denote the element in the $i$th row and $j$th column as $M_n^{i,\,j}$, where $1\leq i \leq n-1$ and $1\leq j\leq n$.
    \begin{equation} \label{complete_original_matrix}
    M_n=\begin{amatrix}{9}
        1 & 0 & -1 & -1 & -1 & \ldots & -1 & -1 & -1 & 1\\
        -1 & 1 & 0 & -1 & -1 & \ldots & -1 & -1 & -1 & 1\\
        0 & -1 & 2 & 0 & -1 & \ldots & -1 & -1 & -1 & 1 \\
        \vdots & & & & & & & & & \vdots \\
        0 & 0 & 0 & 0 & 0 & \ldots & -1 & n-3 & 0 & 1\\
        0 & 0 & 0 & 0 & 0 & \ldots & 0 & -1 & n-2 & 0
    \end{amatrix}  
    \end{equation}
    
    We now apply Gauss-Jordan Elimination to transform $M_n$ into reduced row echelon form. In particular, we seek the value of $T_n=f_3$. We will denote the reduced row echelon form of $M_n$ as $R_n$.

    We denote the submatrix formed by omitting the last row and column of $M_n$ as $M_n[-1]$. Note that $M_{n-1}$ is very similar to $M_n[-1]$. In particular, $M_n[-1]$ differs from $M_{n-1}$ only in the $(n-1)$th column, which is negated. We proceed with the following algorithm.

    We first add each row to the one beneath it, putting $M_n$ into row echelon form. 
    \begin{equation}\label{complete_simplified_matrix}
    M_n=\begin{amatrix}{9}
        1 & 0 & -1 & -1 & -1 & \ldots & -1 & -1 & -1 & 1\\
        0 & 1 & -1 & -2 & -2 & \ldots & -2 & -2 & -2 & 2\\
        0 & 0 & 1 & -2 & -3 & \ldots & -3 & -3 & -3 & 3 \\
        \vdots & & & & & & & & & \vdots \\
        0 & 0 & 0 & 0 & 0 & \ldots & 0 & 1 & -(n-3) & n-2\\
        0 & 0 & 0 & 0 & 0 & \ldots & 0 & 0 & 1 & n-2
    \end{amatrix}   
    \end{equation}
    
    We next apply the rest of the Gauss-Jordan Elimination algorithm for $M_{n-2}$ on $M_n$ to create a new matrix $M_n'$. The $(n-2)$th and $(n-1)$th columns of $M_n'$ will have arbitrary constants, but the rest of $M_n'$ will be in reduced row echelon form. In other words, we solve for $f_2,\,\ldots,\,f_{n-2}$ in terms of $f_{n-1}$ and $f_n$. Recall that $M_n'^{2,\,n-2}=-R_{n-2}^{2,\,n-2}=-T_{n-2}$. To cancel $M_n'^{2,\,n-2}$, we note that $M_n'^{n-2,\,n-2}=1$. Thus, we add $T_{n-2}$ times the $(n-2)$th row to the second row, leaving $f_3$ in terms of only $f_n$. If $M_n^{n-2,\,n}$ were equal to $ n-3 = M_{n-1}^{n-2,\,n-1}$, we would now have that $M_n'^{2,\,n}=T_{n-1}$. However, $M_n^{n-2,\,n}$ is instead $M_{n-1}^{n-2,\,n-1}+1$. Therefore, adding $T_{n-2}$ times the $(n-2)$th row to the 2nd row adds an additional quantity of $T_{n-2}$ to $M_n'^{2,\,n}$, so that we now have $M_n'^{2,\,n}=T_{n-2}+T_{n-1}$.

    We now proceed with the rest of the Gauss-Jordan Elimination algorithm for $M_{n-1}$ on $M_n'$ to create a new matrix $M_n''$, in which only the $(n-1)$th column is not in reduced row echelon form. In other words, we solve for $f_2,\,\ldots,\,f_{n-1}$ in terms of $f_n$. Recall that $M_n''^{2,\,n-1}=-T_{n-1}$. To cancel $M_n''^{2,\,n-1}$, we note that $M_n''^{n-1,\,n-1}=1$. Thus, we add $T_{n-1}$ times the $(n-1)$th row to the second row, putting the second row completely into reduced row echelon form. Note that $M_n''^{n-1,\,n}=\sum_{i=1}^{n-1} M_n^{i,\,n}=\deg^+ v_1=n-2$. Thus, $$M_n''^{2,\,n}=M_n'^{2,\,n}+T_{n-1}(n-2)=T_{n-2}+T_{n-1}+T_{n-1}(n-2)=T_{n-2}+(n-1)T_{n-1}.$$ 
    
    Since the second row is now in reduced row echelon form, we have $$f_3=T_{n-2}+(n-1)T_{n-1}.$$ 
    
    Since $\gcd(f_1,\,f_3)=\gcd(1,\,f_3)=1$, we have $T_n=f_3=T_{n-2}+(n-1)T_{n-1}$.
    
    For $n\leq 3$, \Cref{vertex3} does not hold. However, we can manually solve System \ref{complete_simplified_system} for $n=1,\,2,\,3$ to find the minimal solution $f_i=1$ for all $i$. Thus, we set $T_1=T_2=T_3=1$. For $n\geq 4$, we use the above recursive formula. By the OEIS Sequence A058279 \cite{oeis}, $T_n = \Theta((n-1)!)$. 

    Recall that $f_i$ is the number of times that $v_i$ fires in the first $T$ rounds. Thus, $T\geq f_i$. Since $f_3\geq T_n$ for $n\geq 4$ and $f_1\geq T_n$ for $n\leq 3$, we have that $T\geq T_n$ for all parallel chip-firing games in which at least one vertex fires each round. Thus, $$T_c(K_n)\geq T_n=\Theta((n-1)!) \implies T_c(K_n) = \Omega((n-1)!). \eqno\qed\phantom\qedhere$$
\end{proof}

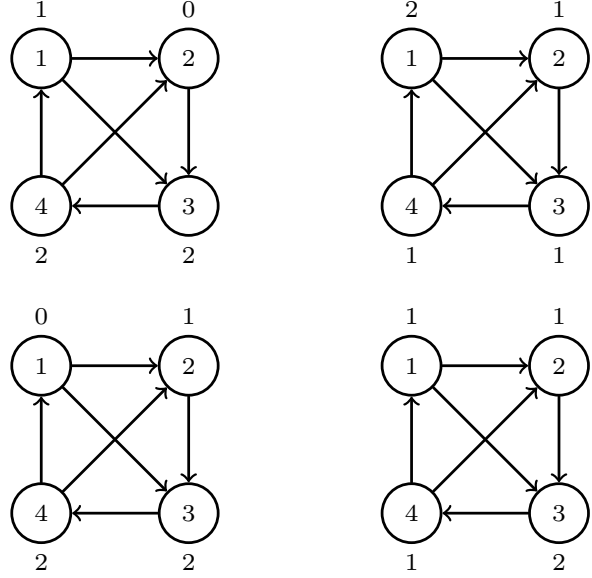
\begin{figure}
    \centering 
    \scalebox{1.3}{\begin{tikzpicture}
      [every node/.style={circle, draw, minimum size=6mm, inner sep=1mm, font=\footnotesize},
       ->, line cap=round, line join=round, thick]
    
      \node (1) at (0, 1.5) {1};
      \node (2) at (1.5, 1.5) {2};
      \node (3) at (1.5, 0) {3};
      \node (4) at (0, 0) {4};
    
      \node[draw=none] at (0, 2) {1};
      \node[draw=none] at (1.5, 2) {0};
      \node[draw=none] at (1.5, -0.5) {2};
      \node[draw=none] at (0, -0.5) {2};
    
      \draw[->] (1) -- (2);
      \draw[->] (2) -- (3);
      \draw[->] (3) -- (4);
      \draw[->] (4) -- (1);
    
      \draw[->] (1) -- (3);
      \draw[->] (4) -- (2);
    
    \end{tikzpicture}}
    \hspace{1.9cm}
    \scalebox{1.3}{\begin{tikzpicture}
      [every node/.style={circle, draw, minimum size=6mm, inner sep=1mm, font=\footnotesize},
       ->, line cap=round, line join=round, thick]
    
      \node (1) at (0, 1.5) {1};
      \node (2) at (1.5, 1.5) {2};
      \node (3) at (1.5, 0) {3};
      \node (4) at (0, 0) {4};
    
      \node[draw=none] at (0, 2) {2};
      \node[draw=none] at (1.5, 2) {1};
      \node[draw=none] at (1.5, -0.5) {1};
      \node[draw=none] at (0, -0.5) {1};
    
      \draw[->] (1) -- (2);
      \draw[->] (2) -- (3);
      \draw[->] (3) -- (4);
      \draw[->] (4) -- (1);
    
      \draw[->] (1) -- (3);
      \draw[->] (4) -- (2);
    
    \end{tikzpicture}}\\
    \scalebox{1.3}{\begin{tikzpicture}
      [every node/.style={circle, draw, minimum size=6mm, inner sep=1mm, font=\footnotesize},
       ->, line cap=round, line join=round, thick]
    
      \node (1) at (0, 1.5) {1};
      \node (2) at (1.5, 1.5) {2};
      \node (3) at (1.5, 0) {3};
      \node (4) at (0, 0) {4};
    
      \node[draw=none] at (0, 2) {0};
      \node[draw=none] at (1.5, 2) {1};
      \node[draw=none] at (1.5, -0.5) {2};
      \node[draw=none] at (0, -0.5) {2};
    
      \draw[->] (1) -- (2);
      \draw[->] (2) -- (3);
      \draw[->] (3) -- (4);
      \draw[->] (4) -- (1);
    
      \draw[->] (1) -- (3);
      \draw[->] (4) -- (2);
    
    \end{tikzpicture}}
    \hspace{1.9cm}
    \scalebox{1.3}{\begin{tikzpicture}
      [every node/.style={circle, draw, minimum size=6mm, inner sep=1mm, font=\footnotesize},
       ->, line cap=round, line join=round, thick]
    
      \node (1) at (0, 1.5) {1};
      \node (2) at (1.5, 1.5) {2};
      \node (3) at (1.5, 0) {3};
      \node (4) at (0, 0) {4};
    
      \node[draw=none] at (0, 2) {1};
      \node[draw=none] at (1.5, 2) {1};
      \node[draw=none] at (1.5, -0.5) {2};
      \node[draw=none] at (0, -0.5) {1};
    
      \draw[->] (1) -- (2);
      \draw[->] (2) -- (3);
      \draw[->] (3) -- (4);
      \draw[->] (4) -- (1);
    
      \draw[->] (1) -- (3);
      \draw[->] (4) -- (2);
    
    \end{tikzpicture}}
    \caption{An example of a parallel chip-firing game on $K_4$ with $T=T_4=4$.}
    \label{example_tournament}
\end{figure}

\begin{example}
    Figure \ref{example_tournament} shows an example of a parallel chip-firing game on $K_4$ with $T=T_4=4$. The vector $(f_1,\,f_2,\,f_3,\,f_4)=(1,\,3,\,4,\,2)$ is indeed a solution to System \ref{complete_original_system} for $K_4$.
\end{example}

\begin{remark}
    Matrix \ref{complete_original_matrix} is extremely similar to the out-degree Laplacian matrix of $K_n$. In fact, the matrix representation of System \ref{complete_original_system} is exactly the out-degree Laplacian. The concept of using the out-degree Laplacian of a directed graph $D$ to analyze the periods of parallel chip-firing games on $D$ was previously mentioned by Prisner \cite{prisner}, though the method of Gauss-Jordan Elimination is novel in the context of chip-firing games.
\end{remark}

\subsection{Complete Bipartite Graphs}
Let $\tau(a,\,b)$ be the maximum possible period of a parallel chip-firing game on an orientation of $K_{a,\,b}$, with $a\leq b$. We show that $\tau(a,\,b)$ is at least asymptotic to $a!$. We conjecture that $\tau(a,\,b)$ is indeed asymptotic to $a!$.

\begin{conjecture}
    \label{conj_bipartite}
    It holds that $\tau(a,\,b)=\Theta(a!)$.
\end{conjecture}

Let the partitions of the complete bipartite graph $K_{a,\,a}$ be $L$ and $R$, where $|L|=|R|=a$. Let $v_i\in L$ if $i\equiv 1\pmod{2}$, and let $v_i\in R$ otherwise. We define the following orientation of $K_{a,a}$ to be a \textit{useful} orientation. 

\begin{itemize}
    \item There exist edges directed from $v_i$ to $v_{i+1}$ for each $1\leq i < n$, and there exists an edge directed from $v_{2a}$ to $v_1$.
    \item There exist edges directed from $v_1$ to $v_{2i}$ for all $1\leq i < a$.
    \item There exist edges directed from $v_{2j+1}$ to all vertices $v_{2k}$ with $1\leq k < j$.
    \item There exist edges directed from $v_{2j}$ to all vertices $v_{2k-1}$ with $1<k<j$.
\end{itemize}

\begin{figure}
    \centering
    \scalebox{1.5}{\begin{tikzpicture}
      [every node/.style={circle, draw, minimum size=6mm, inner sep=1mm, font=\footnotesize},
       ->, line cap=round, line join=round, thick]
    
      \node (1) at (0, 3) {1};
      \node (3) at (0, 2) {3};
      \node (5) at (0, 1) {5};
      \node (7) at (0, 0) {7};
    
      \node (2) at (2, 3) {2};
      \node (4) at (2, 2) {4};
      \node (6) at (2, 1) {6};
      \node (8) at (2, 0) {8};
    
      \draw[->] (1) -- (2);
      \draw[->] (2) -- (3);
      \draw[->] (3) -- (4);
      \draw[->] (4) -- (5);
      \draw[->] (5) -- (6);
      \draw[->] (6) -- (7);
      \draw[->] (7) -- (8);
      \draw[->] (8) -- (1);
    
      \draw[->] (1) -- (4);
      \draw[->] (1) -- (6);
      \draw[->] (6) -- (3);
      \draw[->] (8) -- (3);
      \draw[->] (8) -- (5);
      \draw[->] (5) -- (2);
      \draw[->] (7) -- (2);
      \draw[->] (7) -- (4);
    
    \end{tikzpicture}}
    \caption{A useful orientation of $K_{4,\,4}$.}
    \label{bipartite}
\end{figure}
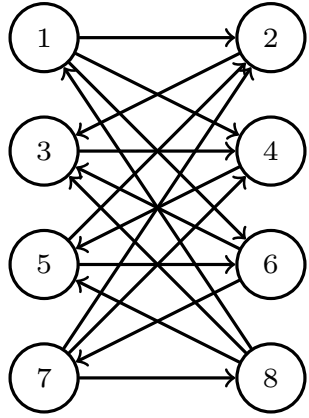

See Figure \ref{bipartite} for a diagram of a useful orientation of $K_{4,\,4}$. Unless otherwise stated, all future mentions of $K_{a,\,a}$ will refer specifically to the useful orientation of $K_{a,\,a}$.

The following lemma shows that vertex 4 fires the most in all parallel chip-firing games on a useful orientation of $K_{a,\,a}$ with $a\geq 3$.

\begin{lemma}
    \label{vertex4}
    For any $a\geq 3$, all parallel chip-firing games on a useful orientation of $K_{a,\,a}$ satisfy $f_4\geq f_i$ for all $1\leq i\leq 2a$.
\end{lemma}

\begin{proof}
    Since $\deg^+ v_4=1$, it must fire once in each period for every chip it receives in that period. Note that $v_4$ receives chips from all vertices $v_{2i-1}$ with $1\leq i \leq a$ except $v_5$, so $f_4\geq f_{2i-1}$ for all $1\leq i \leq a$. Also note that $\deg^+ v_3=1$ and that $v_3$ receives chips from all vertices $v_{2i}$ with $1\leq i\leq a$ and $i\neq 2$. Thus, $f_4\geq f_3\geq f_{2i}$ for all $1\leq i \leq a$ and $i\neq 2$. Finally, $\deg^+ v_2=1$ and $v_2$ receives chips from $v_5$, so $f_4\geq f_3\geq f_2\geq f_5$.
\end{proof}

Using \Cref{vertex4}, we provide a factorial lower bound for the $c$-convergent period of a useful $K_{a,\,a}$ orientation for sufficiently large $c$.

\begin{theorem}
    \label{bipartites}
    For a useful orientation of $K_{a,\, a}$, with $c=a(a-2)+1$, we have $T_c(K_{a, \,a}) = \Omega(a!)$.
\end{theorem}

\begin{proof}
    Since the number of chips on each vertex does not change after a full period, the following equations must hold.
    \begin{equation} \label{bipartite_original_system}
    \begin{aligned}
       &f_{2a} = (a-1)f_1\\
       &f_1 + f_5 + f_7 + \ldots + f_{2a-1} = f_2\\
       &f_1 + f_{2j-1} + f_{2j+3} + \ldots + f_{2a-1} = (j-1)f_{2j} \quad\text{for}\quad 2 \leq j \leq a-2\\
       &f_{2j-2}+f_{2j+2} + f_{2j+4} + \ldots + f_{2a} = (j-1)f_{2j-1} \quad\text{for}\quad 2 \leq j \leq a-1\\
       &f_1 + f_{2a-3} = (a-2)f_{2a-2}\\
       &f_{2a-2}=(a-1)f_{2a-1}\\
       &f_{2a-1} = (a-1)f_{2a}
    \end{aligned}
    \end{equation}
    
    For any solution, scaling gives a new solution, so System \ref{bipartite_original_system} is not independent. We arbitrarily drop the first equation, leaving a system of $2a-1$ independent equations in $2a$ variables. These remaining equations are linearly independent, as can be seen later in the un-augmented form of Matrix \ref{bipartite_simplified_matrix}, whose rank is equal to the number of rows.
    
    By definition, each of $f_1, f_2, \ldots,\, f_{2a} \geq 0.$ By \Cref{lost}, if any $f_i=0$, we have that $T=1$ and that all vertices are forever passive. However, the largest possible $c$ for which all vertices can be forever passive is $\sum_v (\deg^+v-1)=|E|-|V|=a(a-2)<c$. Thus, we further constrain $f_i>0$ for all $i$. Note that for any solution $(f_1, f_2, \ldots ,\,f_{2a})$, it must be true that $(\lambda f_1, \lambda f_2, \ldots,\, \lambda f_{2a})$ is also a solution. To find the minimal solution, we set $f_1=1$. System \ref{bipartite_original_system} has now been reduced to
    \begin{equation} \label{bipartite_simplified_system}
    \begin{aligned}
       &f_2 - f_5 - f_7 - \ldots - f_{2a-1} = 1\\
       &-f_{2j-1} + (j-1)f_{2j} - f_{2j+3} - \ldots - f_{2a-1} = 1 \quad\text{for}\quad 2 \leq j \leq a-2\\
       &-f_{2j-2} + (j-1)f_{2j-1} - f_{2j+2} - f_{2j+4} - \ldots - f_{2a} = 0 \quad\text{for}\quad 2 \leq j \leq a-1\\
       &-f_{2a-3} + (a-2)f_{2a-2} = 1\\
       &-f_{2a-2} + (a-1)f_{2a-1} = 0\\
       &-f_{2a-1} + (a-1)f_{2a} = 0.
    \end{aligned}
    \end{equation}
    
     Note that these equations are linearly independent. Therefore, we can let a well-defined $T_a$ denote the maximal value of $f_i$ satisfying System \ref{bipartite_simplified_system}. By \Cref{vertex4}, for any solution with $a\geq 3$, we have $f_4\geq f_i$ for all $i$. Therefore, $T_a=f_4$ if $a\geq 3$. We claim that, for $a\geq 3$, $$T_a > a\,T_{a-1}.$$
    
    System \ref{bipartite_simplified_system} can be written as the following augmented matrix $M_a=(A|B)$, where $A$ is an $(2a-1)\times (2a-1)$ matrix and $B$ is a $(2a-1)\times 1$ column vector. We denote the element in the $i$th row and $j$th column as $M_a^{i,\,j}$, where $1\leq i \leq 2a-1$ and $1\leq j\leq 2a$.
    \begin{equation}\label{bipartite_original_matrix}
    M_a=\begin{amatrix}{10}
        1 & 0 & 0 & -1 & 0 & -1 & \ldots & 0 & -1 & 0 & 1\\
        -1 & 1 & 0 & 0 & -1 & 0 & \ldots & -1 & 0 & -1 & 0\\
        0 & -1 & 1 & 0 & 0 & -1 & \ldots & 0 & -1 & 0 & 1\\
        0 & 0 & -1 & 2 & 0 & 0 & \ldots & -1 & 0 & -1 & 0\\
        0 & 0 & 0 & -1 & 2 & 0 & \ldots & 0 & -1 & 0 & 1\\
        0 & 0 & 0 & 0 & -1 & 3 & \ldots & -1 & 0 & -1 & 0\\
        \vdots &&&&&&&&&&\vdots \\
        0 & 0 & 0 & 0 & 0 & 0 & \ldots & a-2 & 0 & 0 & 1\\
        0 & 0 & 0 & 0 & 0 & 0 & \ldots & -1 & a-1 & 0 & 0\\
        0 & 0 & 0 & 0 & 0 & 0 & \ldots & 0 & -1 & a-1 & 0 
    \end{amatrix} 
    \end{equation}   
    
    We now apply Gauss-Jordan Elimination to transform $M_a$ into reduced row echelon form. In particular, we seek the value of $T_a=f_4$, which will be the largest of the $f_i$ by \Cref{vertex4}. 

    We denote the submatrix formed by omitting the last two rows and columns of $M_a$ as $M_a[-2]$. Note that $M_{a-1}$ is extremely similar to $M_a[-2]$. In particular, $M_a[-2]$ differs from $M_{a-1}$ only in the $(2a-2)$th column, which is negated. We proceed with the following algorithm. 

    We first add each row to the one beneath it, putting $M_a$ into row echelon form.
    \begin{equation}\label{bipartite_simplified_matrix}
    M_a=\begin{amatrix}{10}
        1 & 0 & 0 & -1 & 0 & -1 & \ldots & 0 & -1 & 0 & 1\\
        0 & 1 & 0 & -1 & -1 & -1 & \ldots & -1 & -1 & -1 & 1\\
        0 & 0 & 1 & -1 & -1 & -2 & \ldots & -1 & -2 & -1 & 2\\
        0 & 0 & 0 & 1 & -1 & -2 & \ldots & -2 & -2 & -2 & 2\\
        0 & 0 & 0 & 0 & 1 & -2 & \ldots & -2 & -3 & -2 & 3\\
        0 & 0 & 0 & 0 & 0 & 1 & \ldots & -3 & -3 & -3 & 3\\
        \vdots &&&&&&&&&&\vdots \\
        0 & 0 & 0 & 0 & 0 & 0 & \ldots & 1 & -(a-2) & -(a-2) & a-1\\
        0 & 0 & 0 & 0 & 0 & 0 & \ldots & 0 & 1 & -(a-2) & a-1\\
        0 & 0 & 0 & 0 & 0 & 0 & \ldots & 0 & 0 & 1 & a-1 
    \end{amatrix}   
    \end{equation}
    
    We next apply the rest of the Gauss-Jordan Elimination algorithm for $M_{a-1}$ on $M_a$ to create a new matrix $M_a'$. The $(2a-2)$th and $(2a-1)$th columns of $M_a'$ will have arbitrary constants, but the rest of $M_a'$ will be in reduced row echelon form. In other words, we solve for $f_2,\,\ldots,\,f_{2a-2}$ in terms of $f_{2a-1}$ and $f_{2a}$. 
    
    Note that the transformation from $M_a$ to $M_a'$ only requires adding positive multiples of rows to other rows, as each row and column of $M_a$ contains only one positive entry. Also, note that the $(2a)$th column of $M_a$ differs from the $(2a-2)$th column in only two ways: they are negations and $|M_a^{2a-3,\,2a}|=|M_a^{2a-3,\,2a-2}|+1$. Therefore, $|M_a'^{i,\,2a}|\geq |M_a'^{i,\,2a-2}|$.
    
    Recall that $M_a'^{3,\,2a-2}=-M_{a-1}^{3,\,2a-2}=-T_{a-1}$. Then we have $M_a'^{3,\,2a}\geq T_{a-1}$. To cancel $M_a'^{3,\,2a-2}$, we note that $M_a'^{2a-2,\,2a-2}=1$. Thus, we add $T_{a-1}$ times the $(2a-2)$th row to the third row, leaving $f_4$ in terms of only $f_{2a}.$ Note that the $2a-2$th row was not altered by the transformation from $M_a$ to $M_a'$, so $M_a'^{2a-2,\,2a}=a-1$. Therefore, we add $(a-1)T_{a-1}$ to $M_a'^{3,\,2a}\geq T_{a-1}$. 

    Note that we now have $M_a'^{3,\,2a-1}=-x$ for some $x$, which we will cancel by adding the $(2a-1)$th row to the third row $x$ times. Importantly, the value of $M_a'^{3,\,2a}$ is not decreased by addition with the $(2a-1)$th row, since $M_a'^{2a-1,\,2a}>0$. In other words, if $M_a''$ is the resultant matrix after the completed Gauss-Jordan Elimination algorithm, $M_a''^{3,\,2a}>M_a'^{3,\,2a}>(a-1)T_{a-1}+T_{a-1}=a\,T_{a-1}.$ Since $f_4=M_a''^{3,\,2a}$ and $\gcd(f_1,\,f_4)=\gcd(1,\,f_4)=1$, we have $$T_a>a\,T_{a-1}.$$ 
    
    For $a\leq 2$, \Cref{vertex4} does not hold. However, we can manually solve System \ref{bipartite_simplified_system} for $a=1,\,2$ to find the minimal solution $f_i=1$ for all $i$. Thus, we set $T_1=T_2=1$. For $a\geq 3$, we have the above inequality.
    
    Recall that $f_i$ is the number of times that $v_i$ fires in the first $T$ rounds. Thus, $T\geq f_i$. Since $f_4\geq T_a$ for $a\geq 3$ and $f_1\geq T_a$ for $a\leq 2$, we have that $T\geq T_a$ for all parallel chip-firing games in which at least one vertex fires each round. Thus, $$T_c(K_{a,\,a})\geq T_a= \Omega(a!). \eqno\qed\phantom\qedhere$$ 
\end{proof}

\begin{remark}
    Matrix \ref{bipartite_original_matrix} is extremely similar to the out-degree Laplacian matrix of $K_{a,\,a}$. As in the case of the complete graph, the matrix representation of \ref{bipartite_original_system} is exactly the out-degree Laplacian.
\end{remark}

Now consider a complete bipartite graph $K_{a,b}$ with partitions $L$ and $R$, such that $|L| = a$ and $|R| = b$ with $a \leq b$. Let $v_i \in L$ for $1 \le i \le a$ and $v_i \in R$ for $a+1 \le i \le a+b$. Define an orientation of $K_{a,b}$ to be \textit{useful} if the following statements hold.
\begin{itemize}
    \item The induced subgraph defined by vertices $v_i$ for $1 \le i \le 2a$ is a useful orientation of $K_{a,a}$.
    \item The induced subgraph defined by vertices $v_i$ for $1 \le i \le 2a$ is the only sink component.
\end{itemize}

\begin{corollary}
    \label{aneqb}
    For the useful orientation of $K_{a,b}$, with $c=(b-a)(a-1)+a(a-2)+1$, we have $T_c(K_{a,\,b}) = \Omega(a!)$.
\end{corollary}
\begin{proof}
   By \Cref{nonsink}, no vertices $v_i$ with $i>2a$ fire in a parallel chip-firing game on $K_{a,\,b}$ so those $b-a$ vertices will have at most $a-1$ chips each. Then the period of the game is simply the period of the induced subgraph $K_{a,a}$, with $c' \ge a(a-2)+1$ chips on it. Thus, $T_c(K_{a,b})=T_{c'}(K_{a,a})=\Omega(a!)$.
\end{proof}

\section{Atomic Firing Sequences}\label{atomic}

Let the \emph{firing sequence} of a vertex $v$ in a parallel chip-firing game on a directed graph be the infinite binary string $F_0(v)F_1(v)\ldots$. Let the \emph{length-$\ell$ firing sequence} of a vertex $v$ be the binary string $F_0(v)F_1(v)\ldots F_{\ell-1}(v)$. Let the \emph{atomic firing sequence} of a vertex $v$ be the length-$T$ firing sequence of $v$. 

We next define some properties of binary strings. Let $s = s_0 s_1 \ldots s_{\ell-1}$ be a length-$\ell$ binary string. A binary string is \emph{periodic with period $k$}, with $1\leq k \leq \lfloor\frac{\ell}{2}\rfloor$, if it can be formed by concatenating copies of a length-$k$ binary string. We first show that no binary string of length $\ell\geq 2$ consisting of only ``0''s can be the atomic firing sequence of a vertex in a parallel chip-firing game on a strongly connected directed graph.

\begin{lemma}
    \label{notall0s}
    For any binary string $s$ of length $\ell\geq 2$ consisting of only ``0''s, no vertex in any parallel chip-firing game on a strongly connected directed graph has $s$ as its atomic firing sequence.
\end{lemma}

\begin{proof}
    If a vertex had $s$ as its atomic firing sequence, it would be forever passive and we would have $T=\ell\geq 2$. However, by \Cref{lost}, each vertex must then fire at least once in each period, so no such vertex exists. 
\end{proof}

We have now established a set of binary strings that cannot be the atomic firing sequence of a vertex in a parallel chip-firing game. We next show that all other binary strings can be atomic firing sequences. First, we prove that all other binary strings of length $\ell \leq 2$ can be atomic firing sequences.

\begin{lemma}
    \label{length1or2}
    For any binary string $s$ of length 1 or 2 excluding $s=00$, there exists a parallel chip-firing game on a directed graph such that some vertex $v$ has $s$ as its atomic firing sequence.
\end{lemma}

\begin{figure}[h]
    \centering
    \begin{tikzpicture}
      [every node/.style={circle, draw, minimum size=6mm, inner sep=2mm},
       ->, >=stealth, line cap=round, line join=round, thick]
    
      \node (1) at (0, 0) {1};
      \node (2) at (2, 0) {2};
      \node (3) at (1, -2) {3};
      \node (4) at (3, -2) {4};

      \node[draw=none] at (0, 0.7) {\small 2};
      \node[draw=none] at (2, 0.7) {\small 2};
      \node[draw=none] at (1, -1.3) {\small 0};
      \node[draw=none] at (3, -1.3) {\small 0};
      
      \draw[->, line width=1.2pt] (1) -- (2);
      \draw[->, line width=1.2pt] (2) -- (3);
      \draw[->, line width=1.2pt] (2) -- (4);
      \draw[->, line width=1.2pt] (3) -- (1);
      \draw[->, line width=1.2pt] (4) -- (1);
    
    \end{tikzpicture}
    \hspace{2cm}
    \begin{tikzpicture}
      [every node/.style={circle, draw, minimum size=6mm, inner sep=2mm},
       ->, >=stealth, line cap=round, line join=round, thick]
    
      \node (1) at (0, 0) {1};
      \node (2) at (2, 0) {2};
      \node (3) at (1, -2) {3};
      \node (4) at (3, -2) {4};

      \node[draw=none] at (0, 0.7) {\small 1};
      \node[draw=none] at (2, 0.7) {\small 1};
      \node[draw=none] at (1, -1.3) {\small 1};
      \node[draw=none] at (3, -1.3) {\small 1};
      
      \draw[->, line width=1.2pt] (1) -- (2);
      \draw[->, line width=1.2pt] (2) -- (3);
      \draw[->, line width=1.2pt] (2) -- (4);
      \draw[->, line width=1.2pt] (3) -- (1);
      \draw[->, line width=1.2pt] (4) -- (1);
    
    \end{tikzpicture}
    \caption{The final parallel chip-firing game from \Cref{length1or2}. Vertex 1 has atomic firing sequence $11$.}
    \label{11}
\end{figure}
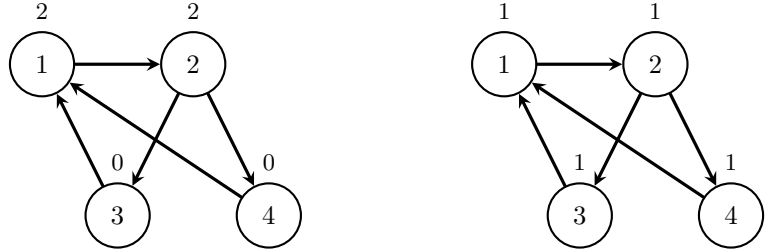

\begin{proof}
    If $s=0$, any vertex in any parallel chip-firing game with no chips will have $s$ as its atomic firing sequence. Similarly, if $s=1$, we can construct the $C_3$ graph, for example, with one chip on each vertex. Then every vertex fires on every round, so every vertex again has $s$ as its atomic firing sequence.

    If $s=01$ or $s=10$, we can construct the $C_4$ graph and place one chip each on $v_0$ and $v_2$. Then $v_0$ and $v_2$ have atomic firing sequence $10$ and $v_1$ and $v_3$ have atomic firing sequence $01$. Finally, if $s=11$, we construct a graph $D$ with 4 vertices and edge set $E=\{(v_1,\,v_2),\,(v_2,\,v_3),\,(v_2,\,v_4),\,(v_3,\,v_1),\,(v_4,\,v_1)\}$ (see \Cref{11}). The parallel chip-firing game on $D$ with $c_0(v_1)=c_0(v_2)=2$ and $c_0(v_3)=c_0(v_4)=0$ has $T=2$, and $v_1$ fires on both rounds. Thus, $v_1$ has atomic firing sequence $s=11$.
\end{proof}

\Cref{length1or2} concerns very short binary strings, for which our main theorem does not hold. However, all longer binary strings that contain at least one 1 can also be the atomic firing sequences of a vertex in a parallel chip-firing game on a directed graph. The following lemma constructs a parallel chip-firing game on a directed cycle in which a vertex has an arbitrary binary string $s$ of length $\ell \geq 3$ as its length-$\ell$ firing sequence.

\begin{lemma}
    \label{non-periodic}
    For any binary string $s$ of length $\ell\geq 3$, there exists a parallel chip-firing game on a directed graph such that a vertex $v$ has $s$ as its length-$\ell$ firing sequence. If $s$ is not periodic, the atomic firing sequence of $v$ is its length-$\ell$ firing sequence.
\end{lemma}

\begin{proof}
    Consider the simple cycle $C_\ell$. Let $c_0(v_i)=s_i$. On each round, all vertices with a chip fire, so the chips shift one vertex along the cycle. Therefore, $c_t(v_0)=s_t$. Since $c_t(v_0)=F_t(v_0)$, we have that $s$ is the length-$\ell$ firing sequence of $v_0$.

    Note that if $s$ is not periodic, the first round $T$ where $c_T(v_i)=c_0(v_i)=s_i$ is $T=\ell$. In this case, the atomic firing sequence of $v_0$ is its length-$\ell$ firing sequence $s$.
\end{proof}

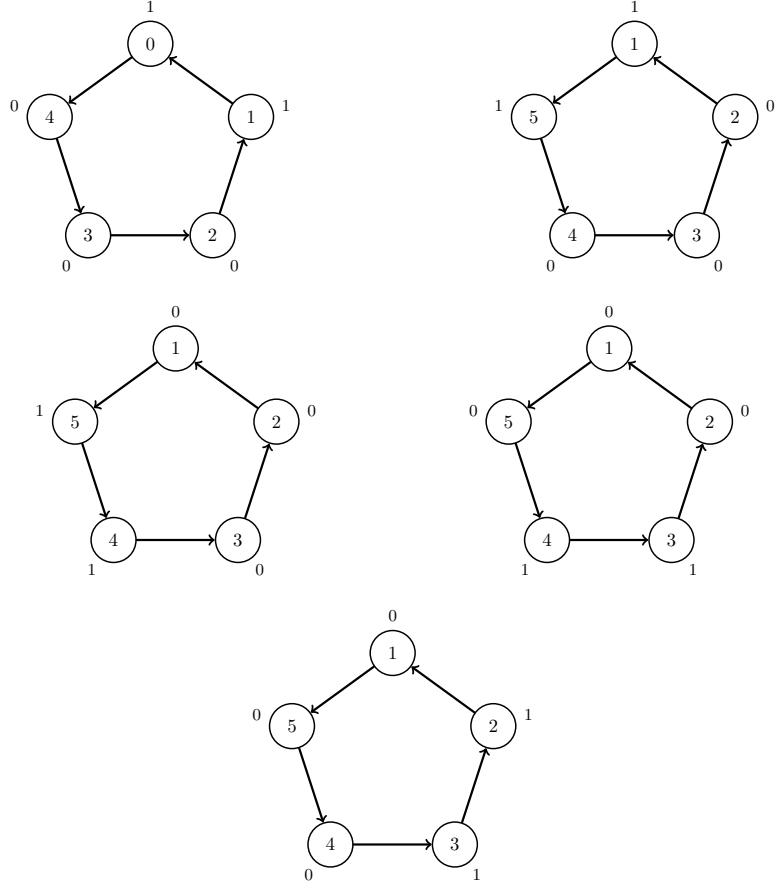
\begin{figure}[h]
    \centering
    \scalebox{0.7}{\begin{tikzpicture}
      [every node/.style={circle, draw, minimum size=6mm, inner sep=2mm},
       ->, line cap=round, line join=round, thick]

      \node (1) at (90:2) {0}; 
      \node (5) at (162:2) {4};
      \node (4) at (234:2) {3}; 
      \node (3) at (306:2) {2}; 
      \node (2) at (18:2) {1};  
    
      \node[draw=none] at (90:2.7) {\small 1};
      \node[draw=none] at (18:2.7) {\small 1};
      \node[draw=none] at (162:2.7) {\small 0};
      \node[draw=none] at (234:2.7) {\small 0};
      \node[draw=none] at (306:2.7) {\small 0};
    
      \draw[<-, line width=1.2pt] (1) -- (2);
      \draw[<-, line width=1.2pt] (2) -- (3);
      \draw[<-, line width=1.2pt] (3) -- (4);
      \draw[<-, line width=1.2pt] (4) -- (5);
      \draw[<-, line width=1.2pt] (5) -- (1);
    
    \end{tikzpicture}}
    \hspace{2cm}
    \scalebox{0.7}{\begin{tikzpicture}
      [every node/.style={circle, draw, minimum size=6mm, inner sep=2mm},
       ->, line cap=round, line join=round, thick]
    
      \node (1) at (90:2) {1}; 
      \node (5) at (162:2) {5}; 
      \node (4) at (234:2) {4}; 
      \node (3) at (306:2) {3}; 
      \node (2) at (18:2) {2};  
    
      \node[draw=none] at (90:2.7) {\small 1};
      \node[draw=none] at (18:2.7) {\small 0};
      \node[draw=none] at (162:2.7) {\small 1};
      \node[draw=none] at (234:2.7) {\small 0};
      \node[draw=none] at (306:2.7) {\small 0};
    
      \draw[<-, line width=1.2pt] (1) -- (2);
      \draw[<-, line width=1.2pt] (2) -- (3);
      \draw[<-, line width=1.2pt] (3) -- (4);
      \draw[<-, line width=1.2pt] (4) -- (5);
      \draw[<-, line width=1.2pt] (5) -- (1);
    
    \end{tikzpicture}}
    \hspace{2cm}
    \scalebox{0.7}{\begin{tikzpicture}
      [every node/.style={circle, draw, minimum size=6mm, inner sep=2mm},
       ->, line cap=round, line join=round, thick]
    
      \node (1) at (90:2) {1};
      \node (5) at (162:2) {5}; 
      \node (4) at (234:2) {4};
      \node (3) at (306:2) {3}; 
      \node (2) at (18:2) {2}; 
    
      \node[draw=none] at (90:2.7) {\small 0};
      \node[draw=none] at (18:2.7) {\small 0};
      \node[draw=none] at (162:2.7) {\small 1};
      \node[draw=none] at (234:2.7) {\small 1};
      \node[draw=none] at (306:2.7) {\small 0};
    
      \draw[<-, line width=1.2pt] (1) -- (2);
      \draw[<-, line width=1.2pt] (2) -- (3);
      \draw[<-, line width=1.2pt] (3) -- (4);
      \draw[<-, line width=1.2pt] (4) -- (5);
      \draw[<-, line width=1.2pt] (5) -- (1);
    
    \end{tikzpicture}
    \hspace{2cm}
    \begin{tikzpicture}
      [every node/.style={circle, draw, minimum size=6mm, inner sep=2mm},
       ->, line cap=round, line join=round, thick]
    
      \node (1) at (90:2) {1};  
      \node (5) at (162:2) {5}; 
      \node (4) at (234:2) {4}; 
      \node (3) at (306:2) {3}; 
      \node (2) at (18:2) {2}; 
    
      \node[draw=none] at (90:2.7) {\small 0};
      \node[draw=none] at (18:2.7) {\small 0};
      \node[draw=none] at (162:2.7) {\small 0};
      \node[draw=none] at (234:2.7) {\small 1};
      \node[draw=none] at (306:2.7) {\small 1};
    
      \draw[<-, line width=1.2pt] (1) -- (2);
      \draw[<-, line width=1.2pt] (2) -- (3);
      \draw[<-, line width=1.2pt] (3) -- (4);
      \draw[<-, line width=1.2pt] (4) -- (5);
      \draw[<-, line width=1.2pt] (5) -- (1);
    
    \end{tikzpicture}}
    \hspace{2cm}
    \scalebox{0.7}{\begin{tikzpicture}
      [every node/.style={circle, draw, minimum size=6mm, inner sep=2mm},
       ->, line cap=round, line join=round, thick]
    
      \node (1) at (90:2) {1}; 
      \node (5) at (162:2) {5};
      \node (4) at (234:2) {4};
      \node (3) at (306:2) {3}; 
      \node (2) at (18:2) {2};  
    
      \node[draw=none] at (90:2.7) {\small 0};
      \node[draw=none] at (18:2.7) {\small 1};
      \node[draw=none] at (162:2.7) {\small 0};
      \node[draw=none] at (234:2.7) {\small 0};
      \node[draw=none] at (306:2.7) {\small 1};
    
      \draw[<-, line width=1.2pt] (1) -- (2);
      \draw[<-, line width=1.2pt] (2) -- (3);
      \draw[<-, line width=1.2pt] (3) -- (4);
      \draw[<-, line width=1.2pt] (4) -- (5);
      \draw[<-, line width=1.2pt] (5) -- (1);
    
    \end{tikzpicture}}
    \caption{An example of a parallel chip-firing game from \Cref{non-periodic}. Here $s=11000$ is not periodic, and $v_0$ has $s$ as its atomic firing sequence.}
    \label{non-periodicexample}
\end{figure}

\begin{example}
    In Figure \ref{non-periodicexample}, $q^0=s=11000=c_0(v_1)c_0(v_2)c_0(v_3)c_0(v_4)c_0(v_5)$. Since $s=11000$ is not periodic, $T=5$.
\end{example}

\Cref{non-periodic} fails to construct a parallel chip-firing game in which a vertex has $s$ as its atomic firing sequence if $s$ is periodic with period $k$, in which case the construction yields a parallel chip-firing game with $T=k<\ell$. To show that a periodic binary string can be the atomic firing sequence of a vertex, we turn to our main theorem.

Let $\ell(s)$ denote the length of a binary string $s$. Let $n(s)$ denote the number of 1s in $s$. Furthermore, let $k(s)$ with $0\leq k(s)\leq \ell(s)-1$ be the maximum value such that $s_{k(s)}=1$. Let $d(s)=\ell(s)-1-k(s)$ be the  number of trailing zeros in $s$. 

\begin{example}
If $s=1010$, then $n(s)=2$ and $k(s)=2$ and $d(s)=1$. If $s=111000$, then $n(s)=3$ and $k(s)=2$ and $d(s)=3$.
\end{example}

In the following proof, we augment the construction from \Cref{non-periodic} to construct a parallel chip-firing game on a directed graph in which a vertex has a periodic atomic firing sequence. To do so, we draw inspiration from Jiang, Scully, and Zhang's \cite{scully} proof of Theorem 4.1 and construct a parallel chip-firing game by attaching many copies of a graph to a single vertex.

\begin{theorem}
    \label{main_atomic}
    For any binary string $s$ that contains 1, there exists a parallel chip-firing game on a directed graph such that $s$ is the firing sequence of some vertex $v$.
\end{theorem}

\begin{proof}
    For brevity, we let $\ell=\ell(s)$ and $n=n(s)$. We also let $k=k(s)$ and $d=d(s)$. If $\ell\leq 2$, we refer to \Cref{length1or2}. 
    
    We create $2n$ copies of the simple cycle $C_\ell$ and label them $C^0,\,C^1,\,\ldots,\,C^{2n-1}$. We denote the $i$th vertex in $C^j$ as $v_i^j$. We assign an initial configuration of chips in each cycle in the same way as in \Cref{non-periodic}: $c_0(v_i^j)=s_i$.  Additionally, we replace all the $v_0^j$ with a single vertex $v$. Thus, there exist edges directed from $v_1^j$ to $v$ and from $v$ to $v_{\ell-1}^j$ for all $j$. Also, $v$ has 0 chips if $s_0=0$, and $v$ has $2n$ chips if $s_0=1$. 

    Since $c_0(v_i^j)$ is independent of $j$, it follows that $c_t(v_i^j)$ is independent of $j$ for all $t$. Thus, collecting all the $v_0^j$ into one vertex does not change $F_t(v_i^j)$ for any $t$ or $v_i^j$. By \Cref{non-periodic}, $v$ has $s$ as its length-$\ell$ firing sequence. We now augment the graph to force $T$ to be $\ell$, regardless of whether $s$ is periodic. 

    We add a new vertex, $u$, with an edge directed from $v$ to $u$. We also add $n$ directed paths of $d$ vertices each. We will call these paths \emph{waterfalls}. If $d=0$, we instead add $n$ single vertices (as if $d=1$). We label the $i$th vertex in the $j$th waterfall $u_i^j$, where $0\leq i\leq d-1$ and $0\leq j\leq n-1$. We add edges directed from $u$ to $u_0^j$ and from $u_{d-1}^j$ to $v$ for all $j$.
    
    If $d=0$, we set $c_0(u)=n$ and $c_0(u_0^j)=0$ for all $j$. If $d\geq 1$, we instead set $c_0(u_{d-1}^j)=1$ and $c_0(u)=c_0(u_i^j)=0$ for all $i\leq d-2$ and all $j$. In both cases, we also place $n$ more chips onto $v$, so that $c_0(v)=s_0(2n)+n$. 

    If $d\geq 1$, on the first round of the game, each vertex $u_{d-1}^j$ fires its 1 chip to $v$. This will give $v$ a total of $2n$ chips, ignoring the behavior of the cycles. Since $\deg^+v=2n+1$, we see that $v$ will only have enough chips to fire when $s_t=1$, so it will still have $s$ as its length-$\ell$ firing sequence. On round $k=\ell-d-1$, vertex $v$ will fire for the $n$th time, so that $c_{k+1}(u)=n=\deg^+ u$. Then $u$ will fire one chip to each $u_0^j$ on round $k+1$. Over the next $d-1$ rounds, these $n$ chips will flow down the waterfalls. Then on round $\ell$, we will have $c_t(u_i^j)=c_0(u_i^j)=1$. We also have $c_\ell(v_i^j)=c_0(v_i^j)$ and $c_\ell(v)=c_0(v)$ by the cyclic construction of the $v_i^j$. Finally, $c_\ell(u)=c_0(u)=0$. Thus, we have $T=\ell$, so $v$ has $s$ as its atomic firing sequence.

    If $d=0$, the game proceeds almost identically. However, $k=\ell-1$, so $c_\ell(u)=c_0(u)=n$. It can be verified that all other vertices have the same number of chips on round $\ell$ as on round 0. Then we again have $T=\ell$, so $v$ has $s$ as its atomic firing sequence.
\end{proof}

\begin{figure}[h]
    \centering
    \scalebox{0.5}{\begin{tikzpicture}
      [every node/.style={circle, draw, minimum size=6mm, inner sep=1mm, font=\footnotesize},
       <-, line cap=round, line join=round, thick]
    
      \node[minimum size=10mm, font=\large] (x) at (2, 0) {$v$};
      \node[minimum size=8mm, font=\large] (y) at (4, 0) {$u$};
      \node (z1) at (3, -1) {$u_0^0$};
      \node (z2) at (5, -1) {$u_0^1$};
    
      \node (a2) at (90:1) {$v_1^0$};
      \node (a3) at (180:1) {$v_2^0$};
      \node (a4) at (270:1) {$v_3^0$};

      \node (b2) at (90:2) {$v_1^1$};
      \node (b3) at (180:2) {$v_2^1$};
      \node (b4) at (270:2) {$v_3^1$};

      \node (c2) at (90:3) {$v_1^2$};
      \node (c3) at (180:3) {$v_2^2$};
      \node (c4) at (270:3) {$v_3^2$};

      \node (d2) at (90:4) {$v_1^3$};
      \node (d3) at (180:4) {$v_2^3$};
      \node (d4) at (270:4) {$v_3^3$};

      \draw[<-, line width=1.2pt] (a4) -- (a3);
      \draw[<-, line width=1.2pt] (a3) -- (a2);
      \draw[<-, line width=1.2pt] (a2) -- (x);
      \draw[<-, line width=1.2pt] (x) -- (a4);

      \draw[<-, line width=1.2pt] (b4) -- (b3);
      \draw[<-, line width=1.2pt] (b3) -- (b2);
      \draw[<-, line width=1.2pt] (b2) -- (x);
      \draw[<-, line width=1.2pt] (x) -- (b4);

      \draw[<-, line width=1.2pt] (c4) -- (c3);
      \draw[<-, line width=1.2pt] (c3) -- (c2);
      \draw[<-, line width=1.2pt] (c2) -- (x);
      \draw[<-, line width=1.2pt] (x) -- (c4);

      \draw[<-, line width=1.2pt] (d4) -- (d3);
      \draw[<-, line width=1.2pt] (d3) -- (d2);
      \draw[<-, line width=1.2pt] (d2) -- (x);
      \draw[<-, line width=1.2pt] (x) -- (d4);

      \draw[->, line width=1.2pt] (x) -- (y);

      \draw[->, line width=1.2pt] (y) -- (z1);
      \draw[->, line width=1.2pt] (y) -- (z2);
      \draw[->, line width=1.2pt] (z1) -- (x);
      \draw[->, line width=1.2pt] (z2) -- (x);

      \node[draw=none] at (0, 4.6) {\small 0};
      \node[draw=none] at (-4.6, 0) {\small 1};
      \node[draw=none] at (0, -4.6) {\small 0};
      \node[draw=none] at (2, 0.7) {\small 6};
      \node[draw=none] at (4, 0.7) {\small 0};
      \node[draw=none] at (3, -1.7) {\small 1};
      \node[draw=none] at (5, -1.7) {\small 1};
    \end{tikzpicture}}
    \scalebox{0.5}{\begin{tikzpicture}
      [every node/.style={circle, draw, minimum size=6mm, inner sep=1mm, font=\footnotesize},
       <-, line cap=round, line join=round, thick]

      \node[minimum size=10mm, font=\large] (x) at (2, 0) {$v$};
      \node[minimum size=8mm, font=\large] (y) at (4, 0) {$u$};
      \node (z1) at (3, -1) {$u_0^0$};
      \node (z2) at (5, -1) {$u_0^1$};

      \node (a2) at (90:1) {$v_1^0$};
      \node (a3) at (180:1) {$v_2^0$};
      \node (a4) at (270:1) {$v_3^0$};

      \node (b2) at (90:2) {$v_1^1$};
      \node (b3) at (180:2) {$v_2^1$};
      \node (b4) at (270:2) {$v_3^1$};

      \node (c2) at (90:3) {$v_1^2$};
      \node (c3) at (180:3) {$v_2^2$};
      \node (c4) at (270:3) {$v_3^2$};

      \node (d2) at (90:4) {$v_1^3$};
      \node (d3) at (180:4) {$v_2^3$};
      \node (d4) at (270:4) {$v_3^3$};

      \draw[<-, line width=1.2pt] (a4) -- (a3);
      \draw[<-, line width=1.2pt] (a3) -- (a2);
      \draw[<-, line width=1.2pt] (a2) -- (x);
      \draw[<-, line width=1.2pt] (x) -- (a4);

      \draw[<-, line width=1.2pt] (b4) -- (b3);
      \draw[<-, line width=1.2pt] (b3) -- (b2);
      \draw[<-, line width=1.2pt] (b2) -- (x);
      \draw[<-, line width=1.2pt] (x) -- (b4);

      \draw[<-, line width=1.2pt] (c4) -- (c3);
      \draw[<-, line width=1.2pt] (c3) -- (c2);
      \draw[<-, line width=1.2pt] (c2) -- (x);
      \draw[<-, line width=1.2pt] (x) -- (c4);
 
      \draw[<-, line width=1.2pt] (d4) -- (d3);
      \draw[<-, line width=1.2pt] (d3) -- (d2);
      \draw[<-, line width=1.2pt] (d2) -- (x);
      \draw[<-, line width=1.2pt] (x) -- (d4);

      \draw[->, line width=1.2pt] (x) -- (y);

      \draw[->, line width=1.2pt] (y) -- (z1);
      \draw[->, line width=1.2pt] (y) -- (z2);
      \draw[->, line width=1.2pt] (z1) -- (x);
      \draw[->, line width=1.2pt] (z2) -- (x);

      \node[draw=none] at (0, 4.6) {\small 1};
      \node[draw=none] at (-4.6, 0) {\small 0};
      \node[draw=none] at (0, -4.6) {\small 1};
      \node[draw=none] at (2, 0.7) {\small 3};
      \node[draw=none] at (4, 0.7) {\small 1};
      \node[draw=none] at (3, -1.7) {\small 0};
      \node[draw=none] at (5, -1.7) {\small 0};
    \end{tikzpicture}}
    \scalebox{0.5}{\begin{tikzpicture}
      [every node/.style={circle, draw, minimum size=6mm, inner sep=1mm, font=\footnotesize},
       <-, line cap=round, line join=round, thick]

      \node[minimum size=10mm, font=\large] (x) at (2, 0) {$v$};
      \node[minimum size=8mm, font=\large] (y) at (4, 0) {$u$};
      \node (z1) at (3, -1) {$u_0^0$};
      \node (z2) at (5, -1) {$u_0^1$};

      \node (a2) at (90:1) {$v_1^0$};
      \node (a3) at (180:1) {$v_2^0$};
      \node (a4) at (270:1) {$v_3^0$};

      \node (b2) at (90:2) {$v_1^1$};
      \node (b3) at (180:2) {$v_2^1$};
      \node (b4) at (270:2) {$v_3^1$};

      \node (c2) at (90:3) {$v_1^2$};
      \node (c3) at (180:3) {$v_2^2$};
      \node (c4) at (270:3) {$v_3^2$};

      \node (d2) at (90:4) {$v_1^3$};
      \node (d3) at (180:4) {$v_2^3$};
      \node (d4) at (270:4) {$v_3^3$};

      \draw[<-, line width=1.2pt] (a4) -- (a3);
      \draw[<-, line width=1.2pt] (a3) -- (a2);
      \draw[<-, line width=1.2pt] (a2) -- (x);
      \draw[<-, line width=1.2pt] (x) -- (a4);

      \draw[<-, line width=1.2pt] (b4) -- (b3);
      \draw[<-, line width=1.2pt] (b3) -- (b2);
      \draw[<-, line width=1.2pt] (b2) -- (x);
      \draw[<-, line width=1.2pt] (x) -- (b4);

      \draw[<-, line width=1.2pt] (c4) -- (c3);
      \draw[<-, line width=1.2pt] (c3) -- (c2);
      \draw[<-, line width=1.2pt] (c2) -- (x);
      \draw[<-, line width=1.2pt] (x) -- (c4);

      \draw[<-, line width=1.2pt] (d4) -- (d3);
      \draw[<-, line width=1.2pt] (d3) -- (d2);
      \draw[<-, line width=1.2pt] (d2) -- (x);
      \draw[<-, line width=1.2pt] (x) -- (d4);

      \draw[->, line width=1.2pt] (x) -- (y);

      \draw[->, line width=1.2pt] (y) -- (z1);
      \draw[->, line width=1.2pt] (y) -- (z2);
      \draw[->, line width=1.2pt] (z1) -- (x);
      \draw[->, line width=1.2pt] (z2) -- (x);

      \node[draw=none] at (0, 4.6) {\small 0};
      \node[draw=none] at (-4.6, 0) {\small 1};
      \node[draw=none] at (0, -4.6) {\small 0};
      \node[draw=none] at (2, 0.7) {\small 7};
      \node[draw=none] at (4, 0.7) {\small 1};
      \node[draw=none] at (3, -1.7) {\small 0};
      \node[draw=none] at (5, -1.7) {\small 0};
    \end{tikzpicture}}
    \scalebox{0.5}{\begin{tikzpicture}
      [every node/.style={circle, draw, minimum size=6mm, inner sep=1mm, font=\footnotesize},
       <-, line cap=round, line join=round, thick]

      \node[minimum size=10mm, font=\large] (x) at (2, 0) {$v$};
      \node[minimum size=8mm, font=\large] (y) at (4, 0) {$u$};
      \node (z1) at (3, -1) {$u_0^0$};
      \node (z2) at (5, -1) {$u_0^1$};

      \node (a2) at (90:1) {$v_1^0$};
      \node (a3) at (180:1) {$v_2^0$};
      \node (a4) at (270:1) {$v_3^0$};

      \node (b2) at (90:2) {$v_1^1$};
      \node (b3) at (180:2) {$v_2^1$};
      \node (b4) at (270:2) {$v_3^1$};

      \node (c2) at (90:3) {$v_1^2$};
      \node (c3) at (180:3) {$v_2^2$};
      \node (c4) at (270:3) {$v_3^2$};

      \node (d2) at (90:4) {$v_1^3$};
      \node (d3) at (180:4) {$v_2^3$};
      \node (d4) at (270:4) {$v_3^3$};

      \draw[<-, line width=1.2pt] (a4) -- (a3);
      \draw[<-, line width=1.2pt] (a3) -- (a2);
      \draw[<-, line width=1.2pt] (a2) -- (x);
      \draw[<-, line width=1.2pt] (x) -- (a4);

      \draw[<-, line width=1.2pt] (b4) -- (b3);
      \draw[<-, line width=1.2pt] (b3) -- (b2);
      \draw[<-, line width=1.2pt] (b2) -- (x);
      \draw[<-, line width=1.2pt] (x) -- (b4);

      \draw[<-, line width=1.2pt] (c4) -- (c3);
      \draw[<-, line width=1.2pt] (c3) -- (c2);
      \draw[<-, line width=1.2pt] (c2) -- (x);
      \draw[<-, line width=1.2pt] (x) -- (c4);

      \draw[<-, line width=1.2pt] (d4) -- (d3);
      \draw[<-, line width=1.2pt] (d3) -- (d2);
      \draw[<-, line width=1.2pt] (d2) -- (x);
      \draw[<-, line width=1.2pt] (x) -- (d4);

      \draw[->, line width=1.2pt] (x) -- (y);

      \draw[->, line width=1.2pt] (y) -- (z1);
      \draw[->, line width=1.2pt] (y) -- (z2);
      \draw[->, line width=1.2pt] (z1) -- (x);
      \draw[->, line width=1.2pt] (z2) -- (x);

      \node[draw=none] at (0, 4.6) {\small 1};
      \node[draw=none] at (-4.6, 0) {\small 0};
      \node[draw=none] at (0, -4.6) {\small 1};
      \node[draw=none] at (2, 0.7) {\small 2};
      \node[draw=none] at (4, 0.7) {\small 2};
      \node[draw=none] at (3, -1.7) {\small 0};
      \node[draw=none] at (5, -1.7) {\small 0};
    \end{tikzpicture}}

    \caption{A parallel chip-firing game in which $v$ has atomic firing sequence $1010$.}
    \label{periodic1}
\end{figure}
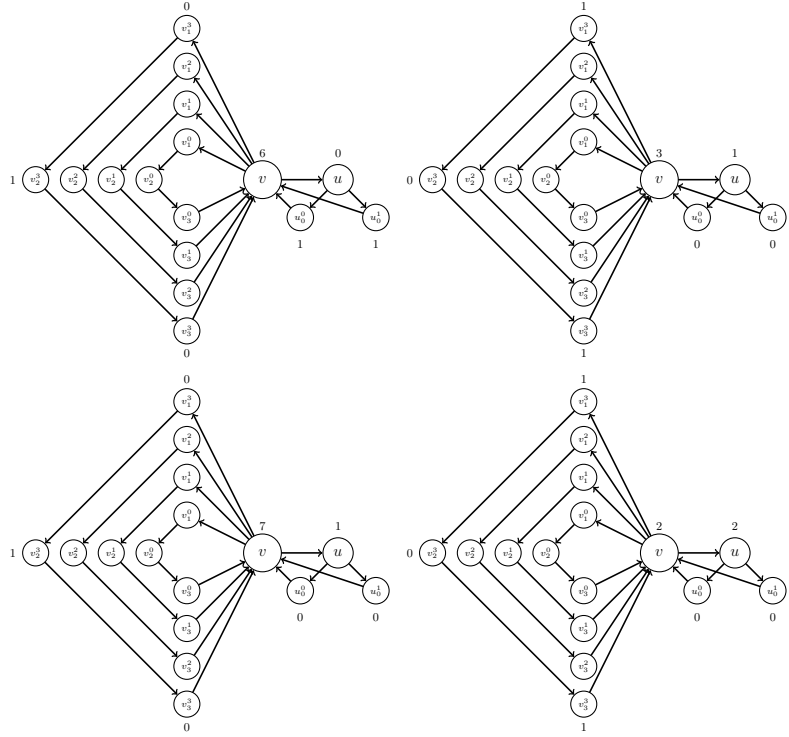

\begin{example}
    \label{atomicexample}
    Figure \ref{periodic1} shows an example construction for $s=1010$. Note that $n=n(s)=2$ and $d=d(s)=1$. Chip numbers written beside the vertices in the outermost cycle apply to the corresponding vertices in each cycle. As shown, $c_0(u_0^0)=c_0(u_0^1)=1$ and $c_0(u)=0$. Also, $c_0(v)=s_0(2n)+n=6$. Since $v$ fires on rounds 0 and 2, we have $k=2$, so $u$ fires on round 3. Then $c_4(u_0^0)=c_4(u_0^1)=1$, so $T=4$. 
\end{example}

\section{Future Directions}\label{futuredirections}

In \Cref{tournaments}, we provided a lower bound for the maximum period of a parallel chip-firing game on an orientation of $K_n$ that is asymptotic to $(n-1)!$. Similarly, in \Cref{bipartites}, we showed that the maximum period of a parallel chip-firing game on an orientation of $K_{a,\,b}$ with $a\leq b$ is at least asymptotic to $a!$. More challenging is the complete characterization of the possible periods of parallel chip-firing games on orientations of $K_n$ or $K_{a,\,b}$. In particular, we conjectured in Conjectures \ref{conj_tourney} and \ref{conj_bipartite} that $\tau(n)$ and $\tau(a,\,b)$ are indeed asymptotic to $(n-1)!$ and $a!$, respectively. Finally, we conjecture that $\lim_{c\to\infty}T_c(D)$ is finite and well-defined for all directed graphs $D$. 

\begin{conjecture}
\label{limctoinfty}
For any directed graph $D$, the quantity $T_c(D)$ is stable for very large $c$. That is, there exists some $T$ and $c_0$ such that $T$ is a possible period of any parallel chip-firing game on $D$ with $c$ chips for all $c\geq c_0$.
\end{conjecture}

For the useful orientations of complete and complete bipartite graphs, we conjecture also that $\lim_{c\to \infty} T_c(D)=T_{c'}(D)$, where $c'$ are the specific values of $c$ defined in Theorems \ref{tournaments} and \ref{bipartites}.

\backmatter

\bmhead{Acknowledgements}

We thank Yunseo Choi for her endless support, for inspiring this project, and for connecting us through our work. 

\bmhead{Data Availability}
No new data were created or analysed. Data sharing is not applicable to this paper.

\begin{appendices}

\section{$T=2$ Games on Undirected Graphs}\label{secA1}

\setcounter{lemma}{0}
\renewcommand{\thelemma}{A\arabic{lemma}}

\begin{lemma}
    \label{talways2}
    For any finite, connected, undirected graph $G=(V,\,E)$ with $|V|\geq 2$, there exists a parallel chip-firing game on $G$ with $T=2$.
\end{lemma}

\begin{proof}
    We set up a parallel chip-firing game on $G$ in the following manner. We arbitrarily choose a vertex $v\in V$ and place $\deg v$ chips onto it. For all other vertices $u$, we use the following process.
    \begin{itemize}
        \item If $d(u,\,v)$ is even, we place $\deg u$ chips onto $u$.
        \item If $d(u,\,v)$ is odd, we  place $\deg u - n$ chips onto $u$, where $n$ is the number of edges $(u,\,w)\in E$ such that $d(u,\,v)\neq d(w,\,v)$. 
    \end{itemize}
    
    On the first round of the game, a vertex $\nu$ fires if and only if $d(\nu,\,v)$ is even. On the second round, all vertices that previously fired have lost chips and thus wait. However, the vertices that previously waited now have exactly enough chips to fire. Thus, all vertices fire exactly once in the first two rounds. By \Cref{constance}, $T=2$.    
\end{proof}




\end{appendices}


\bibliography{sn-bibliography}

\end{document}